\pdfoutput=1
\documentclass[10pt]{amsart} 
\usepackage{mathtools,amssymb,amsthm}

\newtheorem{thm1}{Theorem}[section]
\newtheorem{lem1}[thm1]{Lemma}
\newtheorem{rem1}[thm1]{Remark}

\newtheorem{cor1}[thm1]{Corollary}

\newtheorem{prop1}[thm1]{Proposition}
\newtheorem{ex1}[thm1]{Example}
\newtheorem{not1}[thm1]{Notation}

\begin{document}

\title[]
{Cohen-Macaulayness of associated graded rings of Gorenstein monomial curves}
\author[]{Anargyros Katsabekis}
\address { Department of Mathematics, University of Ioannina, 45110 Ioannina, Greece} \email{katsampekis@uoi.gr}

\begin{abstract} Let $C$ be a Gorenstein noncomplete intersection monomial curve in the 4-dimensional affine space with defining ideal $I(C)$. In this article, we use the minimal generating set of $I(C)$ to give a criterion for determining whether the tangent cone of $C$ is Cohen-Macaulay. We also show that if the tangent cone of $C$ is Cohen-Macaulay, then the minimal number of generators of the ideal $I(C)_{\ast}$  is either $5$ or an even integer of the form $2d+2$, for a suitable integer $d$. Additionally, we provide a family of Gorenstein noncomplete intersection monomial curves $C$ whose tangent cone is Cohen-Macaulay and the minimal number of generators of $I(C)_{\ast}$ is large.
\end{abstract}
\maketitle

\section{Introduction}
It is an important problem to find out whether the associated graded ring of a local ring $(R,\mathfrak{m})$ with respect to its maximal ideal $\mathfrak{m}$ is Cohen-Macaulay. The importance of this problem arises in part from the fact that if the associated graded ring is Cohen-Macaulay, then the computation of the Hilbert function of $R$ can be reduced to the computation of the Hilbert function of an Artin local ring, which is a simple task. In addition, the Cohen-Macaulay property of the associated graded ring guarantees that the Hilbert function is non-decreasing, see \cite{Sta}. This paper addresses the aforementioned problem for local rings associated to Gorenstein noncomplete intersection monomial curves in the 4-dimensional affine space.

Let $\{n_{1},\ldots,n_{r}\}$ be a set of all-different positive integers with greatest common divisor ${\rm gcd}(n_{1},\ldots,n_{r})=1$. Consider the polynomial ring $K[x_{1},\ldots,x_{r}]$ in $r$ variables over a field $K$. Let $C$ be the monomial curve in the $r$-dimensional affine space $A^{r}(K)$ defined parametrically by $x_{1}=t^{n_1},\ldots,x_{r}=t^{n_r}$. The {\em toric ideal} of $C$, denoted by $I(C)$, is the kernel of the $K$-algebra homomorphism $\phi: K[x_{1},\ldots,x_{r}] \rightarrow K[t]$ given by $$\phi(x_{i})=t^{n_i} \ \ \textrm{for all} \ \ 1 \leq i \leq r.$$ It is well known that $I(C)$ is a prime ideal generated by binomials $x_{1}^{u_{1}} \cdots x_{r}^{u_{r}}-x_{1}^{v_{1}} \cdots x_{r}^{v_{r}}$ such that $u_{1}n_{1}+\cdots+u_{r}n_{r}=v_{1}n_{1}+\cdots+v_{r}n_{r}$, see for instance \cite[Lemma 4.1]{Sturmfels95}. Given a polynomial $f \in I(C)$, we let $f_{\ast}$ be the homogeneous summand of $f$ of least degree. We shall denote by $I(C)_{\ast}$ the ideal in $K[x_{1},\ldots,x_{r}]$ generated by the polynomials $f_{\ast}$ for $f \in I(C)$.

Let $\mathfrak{m}=\left< t^{n_1},\ldots,t^{n_r}\right>$ be the maximal ideal of the local ring $R = K[[t^{n_1},\ldots, t^{n_r} ]]$. Then ${\rm gr}_{\mathfrak{m}}(R)=\bigoplus_{i \geq 0} \mathfrak{m}^{i}/\mathfrak{m}^{i+1}$ is the associated graded ring of $R$ and ${\rm gr}_{\mathfrak{m}}(R)$ is isomorphic to the ring $K[x_{1},\ldots,x_{r}]/I(C)_{\ast}$. Thus the Cohen-Macaulayness of ${\rm gr}_{\mathfrak{m}}(R)$ can be studied as the Cohen-Macaulayness of the ring $K[x_1,\ldots, x_r]/I(C)_{\ast}$. Recall that $I(C)_{\ast}$ is the defining ideal of the tangent cone of $C$ at $0$.

We say that the monomial curve $C$ is {\em Gorenstein} if the corresponding local ring $R$ is Gorenstein. For a Gorenstein noncomplete intersection monomial curve $C$ in $A^{4}(K)$, we have from \cite[Theorem 3]{Bresinsky75} that $I(C)$ is minimally generated by the set
$$\{x_3^{a_{13}} x_4^{a_{14}}-x_1^{a_1}, x_{2}^{a_2}- x_{1}^{a_{21}}x_{3}^{a_{23}},  x_3^{a_{3}}-x_{2}^{a_{32}}x_{4}^{a_{34}}, x_{4}^{a_4}-x_{1}^{a_{41}}x_{2}^{a_{42}}, x_{2}^{a_{42}}x_3^{a_{13}}-x_{1}^{a_{21}}x_{4}^{a_{34}}\}$$ where the above binomials are unique up to isomorphism, $a_{ij}>0$ and also $$a_{1} =a_{21}+a_{41}, a_{2}= a_{32}+a_{42}, a_{3}=a_{13}+a_{23}, a_{4} =a_{14}+a_{34}.$$

The paper \cite{AKN} deals with the problem of deciding whether the associated graded ring of a Gorenstein noncomplete intersection local ring $K[[t^{n_1},\ldots, t^{n_4} ]]$ is Cohen-Macaulay by using the minimal generating set of $I(C)$. The authors assume that $n_{1}<n_{2}<n_{3}<n_{4}$ and succeed to provide explicit conditions in many cases, but there are some cases left. More precisely, they do not provide explicit conditions for the following cases: \begin{enumerate} \item Case 1(b) with the restrictions $a_{3}>a_{32}+a_{34}$, $a_{32}<a_{42}$ and $a_{14}>a_{34}$. 
\item Case 2(a)  with the restrictions $a_{24}<a_{34}$ and $a_{13}>a_{23}$. \item Case 2(b) with the restrictions $a_{3}>a_{32}+a_{34}$, $a_{34}<a_{24}$ and $a_{12}>a_{32}$. \item Case 3(b) with the restrictions $a_{23}<a_{43}$ and $a_{14}>a_{24}$.
\end{enumerate}
Our first goal is to provide explicit conditions for the Cohen-Macaulayness of ${\rm gr}_{\mathfrak{m}}(R)$ under the assumption that $n_{1}={\rm min}\{n_{1},\ldots,n_{4}\}$.

Another interesting problem is to estimate the minimal number of generators $\mu(I(C)_{\ast})$ of $I(C)_{\ast}$. In fact J. Herzog and D. Stamate \cite{HerSta} conjectured that, given a monomial curve $C$ in $A^{r}(K)$, ${{\rm wd}(H)+1 \choose 2}$ is an upper bound for $\mu(I(C)_{\ast})$, where ${\rm wd}(H)$ is the width of the semigroup $H=\mathbb{N}\{n_{1},\ldots,n_{r}\}$. Recall that ${\rm wd}(H)$ is the difference between the largest and the smallest generators of $H$. For $r=3$ the conjecture was proved in \cite{LTV}, while it remains open for $r \geq 4$. It is worth noting that T. Shibuta \cite{Shibuta} showed that there is a complete intersection monomial curve $C$ in $A^{4}(K)$ such that ${\rm gr}_{\mathfrak{m}}(R)$ is Cohen-Macaulay and $\mu(I(C)_{\ast})$ is arbitrarily large.

Of particular interest is the case when $C$ is a Gorenstein noncomplete intersection monomial curve in $A^{4}(K)$ and the associated graded ring of $R$ is Cohen-Macaulay. In \cite{AM} F. Arslan and P. Mete showed that $\mu(I(C)_{\ast})=5$ when $I(C)$ is given as in Case 1(a) or 3(a) of \cite{AKN}. In \cite{katsabekis} the author provided families of monomial curves $C$ such that $\mu(I(C)_{\ast})=6$. Our second goal is to provide a formula for $\mu(I(C)_{\ast})$, under the assumption that $C$ is Gorenstein noncomplete intersection and the associated graded ring of $ R$ is Cohen-Macaulay.

In section 2 we assume that $n_{1}={\rm min}\{n_{1},\ldots,n_{4}\}$ and define the number $d$ as the minimum of the set $$\{i \in \mathbb{N}|a_{2}-ia_{32} \leq 0\} \cup \{i \in \mathbb{N}|a_{4}-ia_{34} \leq 0\}.$$ Next we use the minimal generating set of $I(C)$ to provide necessary and sufficient conditions for the Cohen-Macaulayness of ${\rm gr}_{\mathfrak{m}}(R)$ in all possible cases, namely $a_{3} \leq a_{32}+a_{34}$ and $a_{3}>a_{32}+a_{34}$, see Theorems \ref{CM1}, \ref{CM2}. To do so, we first determine a standard basis for $I(C)$ with respect to the negative degree reverse lexicographic term ordering $<$ with $x_{4}>x_{3}>x_{2}>x_{1}$, see Propositions \ref{Basic1}, \ref{Basic2}. Our results apply directly to the Case 1(b) in \cite{AKN} with the restrictions $a_{3}>a_{32}+a_{34}$, $a_{32}<a_{42}$ and $a_{14}>a_{34}$. Similar results can be deduced for the Cases 2(a), 2(b) and 3(b) in \cite{AKN}, see Remarks \ref{basicrem2}, \ref{basicrem3}, \ref{basicrem4}. Additionally, given a Gorenstein noncomplete intersection monomial curve $C$ in $A^{4}(K)$ such that ${\rm gr}_{\mathfrak{m}}(R)$ is Cohen-Macaulay, we show that \[ \mu(I(C)_{\ast})=\begin{cases} 
      5, & \textrm{if} \ a_{3} \leq a_{32}+a_{34} \\
      2d+2, & \textrm{if} \ a_{3}>a_{32}+a_{34}.
   \end{cases}
\]
Finally we provide a family of Gorenstein noncomplete intersection monomial curves $C$ such that ${\rm gr}_{\mathfrak{m}}(R)$ is Cohen-Macaulay and $\mu(I(C)_{\ast})$ is large, see Proposition \ref{infinite}.

\section{Cohen-Macaulay criteria for ${\rm gr}_{\mathfrak{m}}(R)$}

Let $\{n_{1},\ldots,n_{4}\}$  be a set of all-different positive integers with ${\rm gcd}(n_{1},\ldots,n_{4})=1$ and let $C$ be the corresponding monomial curve in $A^{4}(K)$. To prove the main results of this section, namely Theorems \ref{CM1} and \ref{CM2}, we will apply the standard basis algorithm to an appropriate set $S \subset I(C)$. For the definitions of local orderings, normal form, ecart of a polynomial, standard basis and the description of the standard basis algorithm, see \cite{GP}. By using the notation in \cite{GP}, we denote the leading monomial of a polynomial $f$ by ${\rm LM}(f )$, the $S$-polynomial of the polynomials $f$ and $g$ by ${\rm spoly}(f, g)$ and the Mora’s polynomial weak normal form of $f$ with respect to a set of polynomials $S$ by ${\rm NF}(f |S)$. In the sequel, we will make extensive use of the following results.

\begin{lem1} (\cite[Lemma 2.7]{AMS}) \label{AMS} Suppose that $n_{1}={\rm min}\{n_{1},\ldots,n_{4}\}$ and let $S=\{B_{1},\ldots,B_{t}\}$ be a minimal standard basis of the ideal $I(C)$ with respect to a negative degree reverse lexicographical ordering $<$ that makes $x_1$ the lowest variable. Then ${\rm gr}_{\mathfrak{m}}(R)$ is Cohen-Macaulay if and only if $x_{1}$ does not divide ${\rm LM}(B_{i})$ for every $1 \leq i \leq t$.
\end{lem1}

\begin{prop1} (\cite[Proposition 3.1]{Sahin}) \label{Nil} Let $m_{1}$ and $m_{2}$ be monomials and let $B=m_{1}-m_{2}$, $g=m_{1}^{k}-m_{2}^{k}$ where $k \geq 1$ is an integer. If $S$ is a set containing $B$, then ${\rm NF}(g|S)=0$.
\end{prop1}

A {\em binomial} $B=M-N \in I(C)$ is called {\em indispensable} of $I(C)$ if every system of binomial generators of $I(C)$ contains $B$ or $-B$, while a {\em monomial} $M$ is called {\em indispensable} of $I(C)$ if every system of binomial generators of $I(C)$ contains a binomial $B$ such that $M$ is a monomial of $B$. If $C$ is a Gorenstein noncomplete intersection monomial curve, then $I(C)$ is generated by the indispensable binomials, see \cite[Corollary 3.15]{KO}. For the rest of this section we will always assume that $C$ is a Gorenstein noncomplete intersection monomial curve in $A^{4}(K)$. Let $$G:=\{f_{1}=x_3^{a_{13}} x_4^{a_{14}}-x_1^{a_1}, f_{2}=x_{2}^{a_2}- x_{1}^{a_{21}}x_{3}^{a_{23}},  f_{3}=x_3^{a_{3}}-x_{2}^{a_{32}}x_{4}^{a_{34}},$$ $$f_{4}=x_{4}^{a_4}-x_{1}^{a_{41}}x_{2}^{a_{42}}, f_{5}=x_{2}^{a_{42}}x_3^{a_{13}}-x_{1}^{a_{21}}x_{4}^{a_{34}}\}$$ be the unique minimal system of binomial generators of $I(C)$, where $$a_{1} =a_{21}+a_{41}, a_{2}= a_{32}+a_{42}, a_{3}=a_{13}+a_{23}, a_{4} =a_{14}+a_{34}.$$ By \cite[Theorem 4]{Bresinsky75}, $n_{1}=a_{2}a_{4}a_{13}+a_{42}a_{14}a_{23}$, $n_{2}=a_{3}a_{4}a_{21}+a_{41}a_{34}a_{23}$, $n_{3}=a_{1}a_{2}a_{34}+a_{32}a_{21}a_{14}$ and $n_{4}=a_{1}a_{3}a_{42}+a_{13}a_{32}a_{41}$.

First we consider the case that $a_{3} \leq a_{32}+a_{34}$.

\begin{prop1} \label{Basic1} Let $n_{1}={\rm min}\{n_{1},\ldots,n_{4}\}$. Suppose that the following conditions hold.  \begin{enumerate} 
\item $a_{2} \leq a_{21}+a_{23}$.
\item $a_{3} \leq a_{32}+a_{34}$.
\item $a_{4} \leq a_{41}+a_{42}$.
\end{enumerate}
Then $G$ is a minimal standard basis for $I(C)$ with respect to the negative degree reverse lexicographic term ordering $<$ with $x_{4}>x_{3}>x_{2}>x_{1}$.

\end{prop1}
\noindent \textbf{Proof.} Firstly, we show that $a_{42}+a_{13} \leq a_{21}+a_{34}$.  We have that $a_{42}+a_{13}=a_{2}-a_{32}+a_{3}-a_{23}$, $a_{2} \leq a_{21}+a_{23}$ and $a_{3} \leq a_{32}+a_{34}$. Thus $$a_{42}+a_{13} \leq (a_{21}+a_{23})-a_{32}+(a_{32}+a_{34})-a_{23}=a_{21}+a_{34}.$$ Notice that $a_{1}>a_{13}+a_{14}$ since $n_{1}={\rm min}\{n_{1},\ldots,n_{4}\}$.  Here ${\rm LM}(f_{1})=x_3^{a_{13}} x_4^{a_{14}}$, ${\rm LM}(f_{2})=x_{2}^{a_{2}}$, ${\rm LM}(f_{3})=x_{3}^{a_{3}}$, ${\rm LM}(f_{4})=x_{4}^{a_4}$ and ${\rm LM}(f_{5})=x_{2}^{a_{42}}x_3^{a_{13}}$. We have that ${\rm NF}({\rm spoly}(f_{i},f_{j})|G)=0$ as ${\rm LM}(f_{i})$ and ${\rm LM}(f_{j})$ are relatively prime, for $(i,j) \in \{(1,2), (2,3), (2,4), (3,4), (4,5)\}$. Also ${\rm spoly}(f_{1},f_{3})=x_{2}^{a_{32}}x_{4}^{a_{4}}-x_{1}^{a_{1}}x_{3}^{a_{23}}$. Since $a_{4} \leq a_{41}+a_{42}$, we have that $a_{32}+a_{4} \leq a_{41}+a_{2}$ and therefore $a_{32}+a_{4} \leq a_{1}+a_{23}$ since $a_{2} \leq a_{21}+a_{23}$. Thus the leading monomial of ${\rm spoly}(f_{1},f_{3})$ is $x_{2}^{a_{32}}x_{4}^{a_{4}}$, which is divisible only by ${\rm LM}(f_{4})$, and ${\rm ecart}({\rm spoly}(f_{1},f_{3})) \geq {\rm ecart}(f_{4})$. Then ${\rm spoly}(f_{4},{\rm spoly}(f_{1},f_{3}))=x_{1}^{a_{41}}x_{2}^{a_{2}}-x_{1}^{a_{1}}x_{3}^{a_{23}}$ and ${\rm LM}({\rm spoly}(f_{4},{\rm spoly}(f_{1},f_{3})))=x_{1}^{a_{41}}x_{2}^{a_{2}}$. Among the leading monomials of $G$, only ${\rm LM}(f_{2})$ divides $x_{1}^{a_{41}}x_{2}^{a_{2}}$ and ${\rm ecart}({\rm spoly}(f_{4},{\rm spoly}(f_{1},f_{3})))={\rm ecart}(f_2)$. Then ${\rm spoly}(f_{2},{\rm spoly}(f_{4},{\rm spoly}(f_{1},f_{3})))=0$ implying ${\rm NF}({\rm spoly}(f_{1},f_{3})|G)=0$. Also ${\rm spoly}(f_{1},f_{4})=x_{1}^{a_{41}}x_{2}^{a_{42}}x_{3}^{a_{13}}-x_{1}^{a_1}x_{4}^{a_{34}}$ and its leading monomial is $x_{1}^{a_{41}}x_{2}^{a_{42}}x_{3}^{a_{13}}$, which is divided only by ${\rm LM}(f_5)$ and ${\rm ecart}({\rm spoly}(f_{1},f_{4}))={\rm ecart}(f_5)$. Then $${\rm spoly}(f_{5},{\rm spoly}(f_{1},f_{4}))=0$$ implying ${\rm NF}({\rm spoly}(f_{1},f_{4})|G)=0$. Additionally ${\rm spoly}(f_{1},f_{5})=x_{1}^{a_{21}}x_{4}^{a_4}-x_{1}^{a_{1}}x_{2}^{a_{42}}$  and its leading monomial is $x_{1}^{a_{21}}x_{4}^{a_4}$, which is divided only by ${\rm LM}(f_4)$ and $${\rm ecart}({\rm spoly}(f_{1},f_{5}))={\rm ecart}(f_4).$$ Then ${\rm spoly}(f_{4},{\rm spoly}(f_{1},f_{5}))=0$ implying ${\rm NF}({\rm spoly}(f_{1},f_{5})|G)=0$. Furthermore, ${\rm spoly}(f_{2},f_{5})=x_{1}^{a_{21}}x_{2}^{a_{32}}x_{4}^{a_{34}}-x_{1}^{a_{21}}x_{3}^{a_{3}}$  and its leading monomial is $x_{1}^{a_{21}}x_{2}^{a_{32}}x_{4}^{a_{34}}$, which is divided only by ${\rm LM}(f_3)$ and ${\rm ecart}({\rm spoly}(f_{2},f_{5}))={\rm ecart}(f_3)$. Then $${\rm spoly}(f_{3},{\rm spoly}(f_{2},f_{5}))=0$$ implying ${\rm NF}({\rm spoly}(f_{2},f_{5})|G)=0$. Finally ${\rm spoly}(f_{3},f_{5})=x_{1}^{a_{21}}x_{3}^{a_{23}}x_{4}^{a_{34}}-x_{2}^{a_{2}}x_{4}^{a_{34}}$  and its leading monomial is $x_{2}^{a_{2}}x_{4}^{a_{34}}$, which is divided only by ${\rm LM}(f_2)$ and $${\rm ecart}({\rm spoly}(f_{3},f_{5}))={\rm ecart}(f_2).$$ Then ${\rm spoly}(f_{2},{\rm spoly}(f_{3},f_{5}))=0$ implying ${\rm NF}({\rm spoly}(f_{3},f_{5})|G)=0$. \hfill $\square$\\

\begin{thm1} \label{CM1} Let $n_{1}={\rm min}\{n_{1},\ldots,n_{4}\}$. Suppose that $a_{3} \leq a_{32}+a_{34}$. Then ${\rm gr}_{\mathfrak{m}}(R)$ is Cohen-Macaulay if and only if the following conditions hold: \begin{enumerate} 
\item $a_{2} \leq a_{21}+a_{23}$;
\item $a_{4} \leq a_{41}+a_{42}$.
\end{enumerate}
\end{thm1}

\noindent \textbf{Proof.} ($\Leftarrow$) Suppose that (1) and (2) are true. By Proposition \ref{Basic1}, $G$ is a minimal standard basis for $I(C)$ with respect to the negative degree reverse lexicographic term ordering $<$ with $x_{4}>x_{3}>x_{2}>x_{1}$. Since $x_1$ does not divide ${\rm LM}(f_{i})$ for every $1 \leq i \leq 5$, we have from Lemma \ref{AMS} that ${\rm gr}_{\mathfrak{m}}(R)$ is Cohen-Macaulay.\\
 ($\Rightarrow$) Suppose that ${\rm gr}_{\mathfrak{m}}(R)$ is Cohen-Macaulay. The binomials $f_i$, $1 \leq i \leq 5$, are indispensable of $I(C)$, so they belong to a minimal standard basis for $I(C)$ with respect to a negative degree reverse lexicographical ordering $<$ that makes $x_1$ the lowest variable. By Lemma \ref{AMS}, $x_1$ does not divide ${\rm LM}(f_{2})$ and ${\rm LM}(f_{4})$. Thus ${\rm LM}(f_{2})=x_{2}^{a_{2}}$ and ${\rm LM}(f_{4})=x_{4}^{a_4}$. Therefore (1) and (2) are true. \hfill $\square$\\

\begin{ex1} {\rm Consider $n_{1}=141$, $n_{2}=222$, $n_{3}=285$, and $n_{4}=232$. Then $C$ is a Gorenstein noncomplete intersection monomial curve and $I(C)$ is minimally generated by $x_{3}^{3}x_{4}^{3}-x_{1}^{11}, x_{2}^{7}-x_{1}^{9}x_{3}, x_{3}^{4}-x_{2}^{2}x_{4}^{3}, x_{4}^{6}-x_{1}^{2}x_{2}^{5}, x_{2}^{5}x_{3}^{3}-x_{1}^{9}x_{4}^3$. Here $a_{3}=4$, $a_{32}=2$, and $a_{34}=3$ . Note that $a_{2}=7<10=a_{21}+a_{23}$ and $a_{4}=6<7=a_{41}+a_{42}$. By Theorem \ref{CM1}, the ring ${\rm gr}_{\mathfrak{m}}(R)$ is Cohen-Macaulay.}
\end{ex1}

Next we consider the case that $a_{3}>a_{32}+a_{34}$.

\begin{not1} {\rm We shall denote by $d$ the minimum of $$\{i \in \mathbb{N}|a_{2}-ia_{32} \leq 0\} \cup \{i \in \mathbb{N}|a_{4}-ia_{34} \leq 0\}.$$ Notice that $d \geq 2$.}
\end{not1} 
Consider the binomials $$g_{i}=x_{2}^{a_{2}-(i+1)a_{32}}x_{3}^{ia_{3}+a_{13}}-x_{1}^{a_{21}}x_{4}^{(i+1)a_{34}}, 0 \leq i \leq d-2,$$ $$h_{i}=x_{3}^{ia_{3}+a_{13}}x_{4}^{a_{4}-(i+1)a_{34}}-x_{1}^{a_{1}}x_{2}^{ia_{32}}, 0 \leq i \leq d-2,$$ and $$p= \begin{cases} x_{3}^{(d-1)a_{3}+a_{13}}-x_{1}^{a_{21}}x_{2}^{da_{32}-a_{2}}x_{4}^{da_{34}}, \ \ \ \ \ \textrm{if} \ a_{4}-da_{34}> 0\\
x_{3}^{(d-1)a_{3}+a_{13}}-x_{1}^{a_{1}}x_{2}^{(d-1)a_{32}}x_{4}^{da_{34}-a_{4}}, \ \textrm{if} \  a_{4}-da_{34} \leq 0 \end{cases}$$ Note that $g_{0}=f_{5}$ and $h_{0}=f_{1}$. Also the polynomials $g_{i}$, $h_{i}$ and $p$ are in the ideal $I(C)$ by computation, for every $1 \leq i \leq d-2$.

\begin{prop1} \label{Basic2} Let $n_{1}={\rm min}\{n_{1},\ldots,n_{4}\}$. Assume that $a_{3}>a_{32}+a_{34}$ and suppose that the following conditions hold.
\begin{enumerate}
\item $a_{2} \leq a_{21}+a_{23}$.
\item $a_{4} \leq a_{41}+a_{42}$.
\item $a_{42}+a_{13} \leq a_{21}+a_{34}$.
\item $(d-1)(a_{3}-a_{32}-a_{34})+a_{2}-a_{21}-a_{23} \leq 0$.
\item $(d-1)(a_{3}-a_{32}-a_{34})+a_{4}+a_{32}-a_{1}-a_{23} \leq 0$.
\item $d(a_{3}-a_{32}-a_{34})+a_{2}-a_{21}-a_{23} \leq 0$ when $a_{4}-da_{34}>0$.
\item $d(a_{3}-a_{32}-a_{34})+a_{4}+a_{32}-a_{1}-a_{23} \leq 0$ when $a_{4}-da_{34} \leq 0$.
\end{enumerate}
Then $$T=G \cup \{g_{i}| 1 \leq i \leq d-2\} \cup \{h_{i}| 1 \leq i \leq d-2\} \cup \{p\}$$ is a minimal standard basis for $I(C)$ with respect to the negative degree reverse lexicographic term ordering $<$ with $x_{4}>x_{3}>x_{2}>x_{1}$.
\end{prop1}

\noindent \textbf{Proof.} We have that ${\rm LM}(f_{1})=x_3^{a_{13}} x_4^{a_{14}}$, ${\rm LM}(f_{2})=x_{2}^{a_{2}}$, ${\rm LM}(f_{3})=x_{2}^{a_{32}}x_{4}^{a_{34}}$, ${\rm LM}(f_{4})=x_{4}^{a_4}$ and ${\rm LM}(f_{5})=x_{2}^{a_{42}}x_3^{a_{13}}$. For every $1 \leq i \leq d-2$ we have that $i+1 \leq d-1$, thus $(i+1)(a_{3}-a_{32}-a_{34})+a_{2}-a_{21}-a_{23} \leq (d-1)(a_{3}-a_{32}-a_{34})+a_{2}-a_{21}-a_{23}$ since $a_{3}-a_{32}-a_{34}>0$. So $(i+1)(a_{3}-a_{32}-a_{34})+a_{2}-a_{21}-a_{23} \leq 0$ and therefore ${\rm LM}(g_{i})=x_2^{a_{2}-(i+1)a_{32}} x_3^{ia_{3}+a_{13}}$, for every $1 \leq i \leq d-2$. Since $i+1 \leq d-1$ for every $1 \leq i \leq d-2$, we have that $(i+1)(a_{3}-a_{32}-a_{34})+a_{4}+a_{32}-a_{1}-a_{23} \leq (d-1)(a_{3}-a_{32}-a_{34})+a_{4}+a_{32}-a_{1}-a_{23}$ and therefore $(i+1)(a_{3}-a_{32}-a_{34})+a_{4}+a_{32}-a_{1}-a_{23} \leq 0$. Thus ${\rm LM}(h_{i})=x_3^{ia_{3}+a_{13}} x_4^{a_{4}-(i+1)a_{34}}$, for every $1 \leq i \leq d-2$. Additionally ${\rm LM}(p)=x_{3}^{(d-1)a_{3}+a_{13}}$. Since ${\rm LM}(f_{1})$ and  ${\rm LM}(f_{2})$ are relatively prime, we have that ${\rm NF}({\rm spoly}(f_{1},f_{2})|T)=0$. Similarly ${\rm NF}({\rm spoly}(f_{2},f_{4})|T)=0$, ${\rm NF}({\rm spoly}(f_{2},h_{i})|T)=0$, for every $1 \leq i \leq d-2$, and ${\rm NF}({\rm spoly}(f_{4},g_{i})|T)=0$, for every $0 \leq i \leq d-2$.

First we show that ${\rm NF}({\rm spoly}(f_{1},f_{3})|T)=0$. If $a_{4}-2a_{34}>0$, then ${\rm spoly}(f_{1},f_{3})=h_{1}$ and therefore ${\rm NF}({\rm spoly}(f_{1},f_{3})|T)=0$. Suppose that $a_{4}-2a_{34} \leq 0$, so $d=2$. Then ${\rm spoly}(f_{1},f_{3})=p$ and therefore ${\rm NF}({\rm spoly}(f_{1},f_{3})|T)=0$.

Next we show that ${\rm NF}({\rm spoly}(f_{1},f_{4})|T)=0$. We have that ${\rm spoly}(f_{1},f_{4})=x_{1}^{a_{41}}x_{2}^{a_{42}}x_{3}^{a_{13}}-x_{1}^{a_{1}}x_{4}^{a_{34}}$. Then ${\rm LM}( {\rm spoly}(f_{1},f_{4}))=x_{1}^{a_{41}}x_{2}^{a_{42}}x_{3}^{a_{13}}$ and only ${\rm LM}(f_{5})$ divides it. Also ${\rm ecart}({\rm spoly}(f_{1},f_{4}))={\rm ecart}(f_{5})$. Thus ${\rm spoly}(f_{5},{\rm spoly}(f_{1},f_{4}))=0$ and ${\rm NF}({\rm spoly}(f_{1},f_{4})|T)=0$.

We show that ${\rm NF}({\rm spoly}(f_{1},g_{i})|T)=0$, for every $0 \leq i \leq d-2$. We compute ${\rm spoly}(f_{1},g_{i})=x_{1}^{a_{21}}x_{4}^{a_{4}+ia_{34}}-x_{1}^{a_{1}}x_{2}^{a_{2}-(i+1)a_{32}}x_{3}^{ia_{3}}$. Since $a_{32}+a_{34}<a_{3}$ and $a_{4} \leq a_{41}+a_{42}$, we have that ${\rm LM}({\rm spoly}(f_{1},g_{i}))=x_{1}^{a_{21}}x_{4}^{a_{4}+ia_{34}}$ and only ${\rm LM}(f_4)$ divides it. Also ${\rm ecart}({\rm spoly}(f_{1},g_{i})) \geq {\rm ecart}(f_{4})$ and ${\rm spoly}(f_{4},{\rm spoly}(f_{1}, g_{i}))=x_{1}^{a_1}x_{2}^{a_{42}-ia_{32}} \left[(x_{3}^{a_{3}})^{i}-(x_{2}^{a_{32}}x_{4}^{a_{34}})^{i}\right]$. From Proposition \ref{Nil} we have that ${\rm NF}(x_{3}^{ia_{3}}-x_{2}^{ia_{32}}x_{4}^{ia_{34}}|T)$=0 and therefore ${\rm NF}({\rm spoly}(f_{1},g_{i})|T)=0$.

We show that ${\rm NF}({\rm spoly}(f_{1},h_{i})|T)=0$, for every $1 \leq i \leq d-2$. We have that ${\rm spoly}(f_{1},h_{i})=-x_{1}^{a_{1}}\left[(x_{3}^{a_{3}})^{i}-(x_{2}^{a_{32}}x_{4}^{a_{34}})^{i} \right]$. Using Proposition \ref{Nil} we get ${\rm NF}({\rm spoly}(f_{1},h_{i})|T)=0$.

We show that ${\rm NF}({\rm spoly}(f_{1},p)|T)=0$. If $a_{4}-da_{34} \leq 0$, then $${\rm spoly}(f_{1},p)=x_{1}^{a_1}\left[(x_{2}^{a_{32}}x_{4}^{a_{34}})^{d-1}-(x_{3}^{a_{3}})^{d-1} \right].$$  From Proposition \ref{Nil} we have that ${\rm NF}(x_{2}^{(d-1)a_{32}}x_{4}^{(d-1)a_{34}}-x_{3}^{(d-1)a_{3}}|T)$=0 and therefore ${\rm NF}({\rm spoly}(f_{1},p)|T)=0$. Suppose now that $a_{4}-da_{34}>0$. Notice that $da_{32}-a_{2}<a_{32}$ by definition of $d$. Then ${\rm spoly}(f_{1},p)=x_{1}^{a_{21}}x_{2}^{da_{32}-a_{2}}x_{4}^{a_{4}+(d-1)a_{34}}-x_{1}^{a_1}x_{3}^{(d-1)a_{3}}$ and ${\rm LM}({\rm spoly}(f_{1},p))=x_{1}^{a_{21}}x_{2}^{da_{32}-a_{2}}x_{4}^{a_{4}+(d-1)a_{34}}$, since $a_{32}+a_{34}<a_{3}$ and $a_{4} \leq a_{41}+a_{42}$. Only ${\rm LM}(f_4)$ divides ${\rm LM}({\rm spoly}(f_{1},p))$, ${\rm ecart}({\rm spoly}(f_{1},p))> {\rm ecart}(f_4)$ and ${\rm spoly}(f_4,{\rm spoly}(f_1,p))=-x_{1}^{a_1}\left[(x_{2}^{a_{32}}x_{4}^{a_{34}})^{d-1}-(x_{3}^{a_{3}})^{d-1} \right].$ Using Proposition \ref{Nil} we get ${\rm NF}({\rm spoly}(f_{1},p)|T)=0$.

We show that ${\rm NF}({\rm spoly}(f_{2},f_{3})|T)=0$. We have that ${\rm spoly}(f_{2},f_{3})=x_{2}^{a_{42}}x_{3}^{a_{3}}-x_{1}^{a_{21}}x_{3}^{a_{23}}x_{4}^{a_{34}}$, ${\rm LM}({\rm spoly}(f_{2},f_{3}))=x_{2}^{a_{42}}x_{3}^{a_{3}}$ and only ${\rm LM}(f_5)$ divides it. Furthermore ${\rm ecart}({\rm spoly}(f_{2},f_{3}))={\rm ecart}(f_5)$, ${\rm spoly}(f_5,{\rm spoly}(f_{2},f_{3}))=0$ and ${\rm NF}({\rm spoly}(f_{2},f_{3})|T)=0$.

We show that ${\rm NF}({\rm spoly}(f_{2},g_{i})|T)=0$, for every $0 \leq i \leq d-2$. We have that ${\rm spoly}(f_{2},g_{i})=x_{1}^{a_{21}}\left[(x_{2}^{a_{32}}x_{4}^{a_{34}})^{i+1}-(x_{3}^{a_{3}})^{i+1} \right]$. From Proposition \ref{Nil} we have that ${\rm NF}(x_{2}^{(i+1)a_{32}}x_{4}^{(i+1)a_{34}}-x_{3}^{(i+1)a_{3}}|T)$=0 and therefore ${\rm NF}({\rm spoly}(f_{2},g_{i})|T)=0$.

We show that ${\rm NF}({\rm spoly}(f_{3},f_{4})|T)=0$. We have that ${\rm spoly}(f_{3},f_{4})=x_{1}^{a_{41}}x_{2}^{a_2}-x_{3}^{a_3}x_{4}^{a_{14}}$. We distinguish the following cases: \begin{enumerate} \item ${\rm LM}({\rm spoly}(f_{3},f_{4}))=x_{1}^{a_{41}}x_{2}^{a_2}$. Then only ${\rm LM}(f_2)$ divides ${\rm LM}({\rm spoly}(f_{3},f_{4}))$ and ${\rm ecart}({\rm spoly}(f_{3},f_{4}))<{\rm ecart}(f_2)$, since $a_{13}+a_{14}<a_{1}$. Also $${\rm spoly}(f_{2},{\rm spoly}(f_{3},f_{4}))=x_{3}^{a_3}x_{4}^{a_{14}}-x_{1}^{a_1}x_{3}^{a_{23}}.$$ We have that ${\rm LM}({\rm spoly}(f_{2},{\rm spoly}(f_{3},f_{4})))=x_{3}^{a_3}x_{4}^{a_{14}}$ and only ${\rm LM}(f_1)$ divides it. Moreover ${\rm ecart}({\rm spoly}(f_{2},{\rm spoly}(f_{3},f_{4})))={\rm ecart}(f_1)$, $${\rm spoly}(f_1,{\rm spoly}(f_{2},{\rm spoly}(f_{3},f_{4})))=0$$ and ${\rm NF}({\rm spoly}(f_{3},f_{4})|T)=0$. \item ${\rm LM}({\rm spoly}(f_{3},f_{4}))=x_{3}^{a_3}x_{4}^{a_{14}}$. Then only ${\rm LM}(f_1)$ divides ${\rm LM}({\rm spoly}(f_{3},f_{4}))$ and ${\rm ecart}({\rm spoly}(f_{3},f_{4})) \leq {\rm ecart}(f_1)$, since $a_{2} \leq a_{21}+a_{23}$. Also $${\rm spoly}(f_{1},{\rm spoly}(f_{3},f_{4}))=x_{1}^{a_{41}}x_{2}^{a_{2}}-x_{1}^{a_1}x_{3}^{a_{23}}.$$ We have that ${\rm LM}({\rm spoly}(f_{1},{\rm spoly}(f_{3},f_{4})))=x_{1}^{a_{41}}x_{2}^{a_{2}}$ and only ${\rm LM}(f_2)$ divides it. Moreover ${\rm ecart}({\rm spoly}(f_{1},{\rm spoly}(f_{3},f_{4})))={\rm ecart}(f_2)$, $${\rm spoly}(f_2,{\rm spoly}(f_{1},{\rm spoly}(f_{3},f_{4})))=0$$ and ${\rm NF}({\rm spoly}(f_{3},f_{4})|T)=0$.
\end{enumerate}

We show that ${\rm NF}({\rm spoly}(f_{3},g_{i})|T))=0$, for every $0 \leq i \leq d-3$. Then $a_{2}-(i+1)a_{32}>a_{32}$ by definition of $d$, since $i+1<d$, and therefore ${\rm spoly}(f_{3},g_{i})=x_{1}^{a_{21}}x_{4}^{(i+2)a_{34}}-x_{2}^{a_{2}-(i+2)a_{32}}x_{3}^{(i+1)a_{3}+a_{13}}$. Since ${\rm spoly}(f_{3},g_{i})=-g_{i+1}$, we get $${\rm NF}({\rm spoly}(f_{3},g_{i})|T)=0.$$

We show that ${\rm NF}({\rm spoly}(f_{3},g_{d-2})|T)=0$. Suppose first that $a_{4}-da_{34}>0$, then $a_{2}-da_{32} \leq 0$. We have that ${\rm spoly}(f_{3},g_{d-2})=-p$ and therefore ${\rm NF}({\rm spoly}(f_{3},g_{d-2})|T)=0$. Suppose now that $a_{4}-da_{34}\leq 0$. Assume that $a_{2}-da_{32} \leq 0$. Then ${\rm spoly}(f_{3},g_{d-2})=x_{1}^{a_{21}}x_{2}^{da_{32}-a_{2}}x_{4}^{da_{34}}-x_{3}^{(d-1)a_{3}+a_{13}}$ and ${\rm LM}({\rm spoly}(f_{3},g_{d-2}))=x_{3}^{(d-1)a_{3}+a_{13}}$, since $d(a_{3}-a_{32}-a_{34})+a_{2}-a_{21}-a_{23} \leq 0$. Only ${\rm LM}(p)$ divides ${\rm LM}({\rm spoly}(f_{3},g_{d-2}))$ and ${\rm ecart}({\rm spoly}(f_{3},g_{d-2})) \leq {\rm ecart}(p)$, since $a_{4} \leq a_{41}+a_{42}$. We have that $${\rm spoly}(p,{\rm spoly}(f_{3},g_{d-2}))=x_{1}^{a_1}x_{2}^{(d-1)a_{32}}x_{4}^{da_{34}-a_{4}}-x_{1}^{a_{21}}x_{2}^{(d-1)a_{32}-a_{42}}x_{4}^{da_{34}}.$$ Then ${\rm LM}({\rm spoly}(p,{\rm spoly}(f_{3},g_{d-2})))=x_{1}^{a_{21}}x_{2}^{(d-1)a_{32}-a_{42}}x_{4}^{da_{34}}$ and only ${\rm LM}(f_4)$ divides it. Also ${\rm ecart}({\rm spoly}(p,{\rm spoly}(f_{3},g_{d-2})))={\rm ecart}(f_4)$, $${\rm spoly}(f_{4},{\rm spoly}(p,{\rm spoly}(f_{3},g_{d-2})))=0$$ and ${\rm NF}({\rm spoly}(f_{3},g_{d-2})|T)=0$. Assume that $a_{2}-da_{32}>0$. Then ${\rm spoly}(f_{3},g_{d-2})=x_{1}^{a_{21}}x_{4}^{da_{34}}-x_{2}^{a_{2}-da_{32}}x_{3}^{(d-1)a_{3}+a_{13}}$. We distinguish the following cases: \begin{enumerate}
\item ${\rm LM}({\rm spoly}(f_{3},g_{d-2}))=x_{2}^{a_{2}-da_{32}}x_{3}^{(d-1)a_{3}+a_{13}}$. Only ${\rm LM}(p)$ divides the monomial ${\rm LM}({\rm spoly}(f_{3},g_{d-2}))$ and ${\rm ecart}({\rm spoly}(f_{3},g_{d-2})) \leq {\rm ecart}(p)$, since $a_{4} \leq a_{41}+a_{42}$. Then $${\rm spoly}(p,{\rm spoly}(f_{3},g_{d-2}))=x_{1}^{a_{21}}x_{4}^{da_{34}}-x_{1}^{a_1}x_{2}^{a_{42}}x_{4}^{da_{34}-a_{4}}$$ and ${\rm LM}({\rm spoly}(p,{\rm spoly}(f_{3},g_{d-2})))=x_{1}^{a_{21}}x_{4}^{da_{34}}$, since $d(a_{3}-a_{32}-a_{34})+a_{2}-a_{21}-a_{23}$. Only ${\rm LM}(f_4)$ divides ${\rm LM}({\rm spoly}(p,{\rm spoly}(f_{3},g_{d-2})))$ and ${\rm ecart}({\rm spoly}(p,{\rm spoly}(f_{3},g_{d-2})))={\rm ecart}(f_4)$. Thus $${\rm spoly}(f_{4},{\rm spoly}(p,{\rm spoly}(f_{3},g_{d-2})))=0$$ and ${\rm NF}({\rm spoly}(f_{3},g_{d-2})|T)=0$. \item ${\rm LM}({\rm spoly}(f_{3},g_{d-2}))=x_{1}^{a_{21}}x_{4}^{da_{34}}$. Only ${\rm LM}(f_4)$ divides ${\rm LM}({\rm spoly}(f_{3},g_{d-2}))$. We have that ${\rm spoly}(f_{4},{\rm spoly}(f_{3},g_{d-2}))=x_{2}^{a_{2}-da_{32}}x_{3}^{(d-1)a_{3}+a_{13}}-x_{1}^{a_1}x_{2}^{a_{42}}x_{4}^{da_{34}-a_{4}}$ and ${\rm LM}({\rm spoly}(f_{4},{\rm spoly}(f_{3},g_{d-2})))=x_{2}^{a_{2}-da_{32}}x_{3}^{(d-1)a_{3}+a_{13}}$. Only ${\rm LM}(p)$ divides ${\rm LM}({\rm spoly}(f_{4},{\rm spoly}(f_{3},g_{d-2})))$ and ${\rm ecart}({\rm spoly}(f_{4},{\rm spoly}(f_{3},g_{d-2})))={\rm ecart}(p)$. Thus ${\rm spoly}(p,{\rm spoly}(f_{4},{\rm spoly}(f_{3},g_{d-2})))=0$ and $${\rm NF}({\rm spoly}(f_{3},g_{d-2})|T)=0.$$

\end{enumerate}

We show that ${\rm NF}({\rm spoly}(f_{3},h_{i})|T)=0$, for every $1 \leq i \leq d-3$. Then $a_{4}-(i+1)a_{34}>a_{34}$ by definition of $d$, since $i+1<d$, and therefore ${\rm spoly}(f_{3},h_{i})=x_{1}^{a_{1}}x_{2}^{(i+1)a_{32}}-x_{3}^{(i+1)a_{3}+a_{13}}x_{4}^{a_{4}-(i+2)a_{34}}$. Since ${\rm spoly}(f_{3},h_{i})=-h_{i+1}$, we get $${\rm NF}({\rm spoly}(f_{3},h_{i})|T)=0.$$

We show that ${\rm NF}({\rm spoly}(f_{3},h_{d-2})|T)=0$. Suppose first that $a_{4}-da_{34}>0$, so $a_{2}-da_{32} \leq 0$ by definition of $d$. Then ${\rm spoly}(f_{3},h_{d-2})=x_{1}^{a_1}x_{2}^{(d-1)a_{32}}-x_{3}^{(d-1)a_{3}+a_{13}}x_{4}^{a_{4}-da_{34}}$. Then ${\rm LM}({\rm spoly}(f_{3},h_{d-2}))=x_{3}^{(d-1)a_{3}+a_{13}}x_{4}^{a_{4}-da_{34}}$, since $d(a_{3}-a_{32}-a_{34})+a_{2}-a_{21}-a_{23} \leq 0$ and $a_{4} \leq a_{41}+a_{42}$. Only ${\rm LM}(p)$ divides ${\rm LM}({\rm spoly}(f_{3},h_{d-2}))$, ${\rm ecart}({\rm spoly}(f_{3},h_{d-2})) \geq {\rm ecart}(p)$ and $${\rm spoly}(p,{\rm spoly}(f_{3},h_{d-2}))=x_{1}^{a_1}x_{2}^{(d-1)a_{32}}-x_{1}^{a_{21}}x_{2}^{da_{32}-a_{2}}x_{4}^{a_4}.$$ Notice that $da_{32}-a_{2}<a_{32}$ by definition of $d$. Then ${\rm LM}({\rm spoly}(p,{\rm spoly}(f_{3},h_{d-2})))=x_{1}^{a_{21}}x_{2}^{da_{32}-a_{2}}x_{4}^{a_4}$ and only ${\rm LM}(f_{4})$ divides it. Also ${\rm ecart}({\rm spoly}(p,{\rm spoly}(f_{3},h_{d-2})))={\rm ecart}(f_{4})$, ${\rm spoly}(p,{\rm spoly}(f_{3},{\rm spoly}(f_{3},h_{d-2})))=0$ and ${\rm NF}({\rm spoly}(f_{3},h_{d-2})|T)=0$. Suppose now that $a_{4}-da_{34} \leq 0$. Then ${\rm spoly}(f_{3},h_{d-2})=-p$ and therefore ${\rm NF}({\rm spoly}(f_{3},h_{d-2})|T)=0$.

We show that ${\rm NF}({\rm spoly}(f_{4},h_{i})|T)=0$, for every $1 \leq i \leq d-2$.  We have that ${\rm spoly}(f_{4},h_{i})=x_{1}^{a_1}x_{2}^{ia_{32}}x_{4}^{(i+1)a_{34}}-x_{1}^{a_{41}}x_{2}^{a_{42}}x_{3}^{ia_{3}+a_{13}}$ and ${\rm LM}({\rm spoly}(f_{4},h_{i}))=x_{1}^{a_{41}}x_{2}^{a_{42}}x_{3}^{ia_{3}+a_{13}}$, since $i(a_{3}-a_{32}-a_{34})+a_{2}-a_{21}-a_{23} \leq 0$. The elements of the set $Q=\{{\rm LM}(g_{j})| 0 \leq j \leq i\}$ are the divisors of ${\rm LM}({\rm spoly}(f_{4},h_{i}))$. Clearly ${\rm ecart}(g_{i})$ is minimal among the ecarts of the elements of $Q$. Also ${\rm ecart}({\rm spoly}(f_{4},h_{i}))= {\rm ecart}(g_{i})$, ${\rm spoly}(g_{i},{\rm spoly}(f_{4},h_{i}))=0$ and ${\rm NF}({\rm spoly}(f_{4},h_{i})|T)=0$.

We show that ${\rm NF}({\rm spoly}(g_{i},g_{j})|T)=0$, for every $i,j \in \{0,\ldots,d-2\}$ with $i<j$. We have that ${\rm spoly}(g_{i},g_{j})=x_{1}^{a_{21}}x_{4}^{(i+1)a_{34}}\left[(x_{2}^{a_{32}}x_{4}^{a_{34}})^{j-i}-(x_{3}^{a_{3}})^{j-i} \right]$. From Proposition \ref{Nil} we have that ${\rm NF}(x_{2}^{(j-i)a_{32}}x_{4}^{(j-i)a_{34}}-x_{3}^{(j-i)a_{3}}|T)=0$ and also ${\rm NF}({\rm spoly}(g_{i},g_{j})|T)=0$.

We show that ${\rm NF}({\rm spoly}(g_{i},h_{j})|T)=0$, for every $1 \leq i \leq d-2$ and $1 \leq j \leq d-2$. We distinguish the following cases: \begin{enumerate} \item $i<j$. Then ${\rm spoly}(g_{i},h_{j})=x_{1}^{a_1}x_{2}^{a_{2}+(j-i-1)a_{32}}-x_{1}^{a_{21}}x_{3}^{(j-i)a_{3}}x_{4}^{a_{4}-(j-i)a_{34}}$. There are two cases: (a) ${\rm LM}({\rm spoly}(g_{i},h_{j}))=x_{1}^{a_1}x_{2}^{a_{2}+(j-i-1)a_{32}}$. Only ${\rm LM}(f_2)$ divides ${\rm LM}({\rm spoly}(g_{i},h_{j}))$ and $${\rm spoly}(f_{2},{\rm spoly}(g_{i},h_{j}))=x_{1}^{a_{21}}x_{3}^{(j-i)a_{3}}x_{4}^{a_{4}-(j-i)a_{34}}-x_{1}^{a_{1}+a_{21}}x_{2}^{(j-i-1)a_{32}}x_{3}^{a_{23}}.$$ Notice that $a_{4}-(j-i)a_{34}<a_{14}$ for every $j>i+1$. Then $${\rm LM}({\rm spoly}(f_{2},{\rm spoly}(g_{i},h_{j})))=x_{1}^{a_{21}}x_{3}^{(j-i)a_{3}}x_{4}^{a_{4}+(i-j)a_{34}},$$ since $(j-i)(a_{3}-a_{32}-a_{34})+a_{4}+a_{32}-a_{1}-a_{23} \leq 0$, and only ${\rm LM}(h_{j-i-1})$ divides it. Additionally ${\rm ecart}({\rm spoly}(f_{2},{\rm spoly}(g_{i},h_{j}))) ={\rm ecart}(h_{j-i-1})$, $${\rm spoly}(h_{j-i-1},{\rm spoly}(f_{2},{\rm spoly}(g_{i},h_{j})))=0$$ and ${\rm NF}({\rm spoly}(g_{i},h_{j})|T)=0$.\\ (b) ${\rm LM}({\rm spoly}(g_{i},h_{j}))=x_{1}^{a_{21}}x_{3}^{(j-i)a_{3}}x_{4}^{a_{4}-(j-i)a_{34}}$. Only ${\rm LM}(h_{j-i-1})$ divides ${\rm LM}({\rm spoly}(g_{i},h_{j}))$ and ${\rm spoly}(h_{j-i-1},{\rm spoly}(g_{i},h_{j}))=x_{1}^{a_{1}}x_{2}^{a_{2}+(j-i-1)a_{32}}-x_{1}^{a_{1}+a_{21}}x_{2}^{(j-i-1)a_{32}}x_{3}^{a_{23}}$. Then $${\rm LM}({\rm spoly}(h_{j-i-1},{\rm spoly}(g_{i},h_{j})))=x_{1}^{a_{1}}x_{2}^{a_{2}+(j-i-1)a_{32}},$$ since $a_{2} \leq a_{21}+a_{23}$, and only ${\rm LM}(f_2)$ divides it. Also $${\rm ecart}({\rm spoly}(h_{j-i-1},{\rm spoly}(g_{i},h_{j})))={\rm ecart}(f_2),$$ ${\rm spoly}(f_{2},{\rm spoly}(h_{j-i-1},{\rm spoly}(g_{i},h_{j})))=0$ and ${\rm NF}({\rm spoly}(g_{i},h_{j})|T)=0$.

\item $i=j$. Then ${\rm spoly}(g_{i},h_{i})=x_{1}^{a_1}x_{2}^{a_{42}}-x_{1}^{a_{21}}x_{4}^{a_{4}}$, ${\rm LM}({\rm spoly}(g_{i},h_{i}))=x_{1}^{a_{21}}x_{4}^{a_{4}}$ and only ${\rm LM}(f_4)$ divides it. Also ${\rm ecart}({\rm spoly}(g_i,h_i))={\rm ecart}(f_4)$, $${\rm spoly}(f_{4},{\rm spoly}(g_i,h_i))=0$$ and ${\rm NF}({\rm spoly}(g_{i},h_{i})|T)=0$. 

\item $i>j$. Then ${\rm spoly}(g_{i},h_{j})=x_{1}^{a_1}x_{2}^{a_{2}-(i-j+1)a_{32}}x_{3}^{(i-j)a_{3}}-x_{1}^{a_{21}}x_{4}^{a_{4}+(i-j)a_{34}}$ and ${\rm LM}({\rm spoly}(g_{i},h_{j}))=x_{1}^{a_{21}}x_{4}^{a_{4}+(i-j)a_{34}}$, since $a_{4} \leq a_{41}+a_{42}$ and $a_{32}+a_{34}<a_{3}$. Only ${\rm LM}(f_4)$ divides ${\rm LM}({\rm spoly}(g_{i},h_{j}))$ and ${\rm ecart}({\rm spoly}(g_{i},h_{j}))>{\rm ecart}(f_4)$. Also $${\rm spoly}(f_{4},{\rm spoly}(g_{i},h_{j}))=x_{1}^{a_{1}}x_{2}^{a_{2}-(i-j+1)a_{32}}\left[(x_{3}^{a_{3}})^{i-j}-(x_{2}^{a_{32}}x_{4}^{a_{34}})^{i-j}\right].$$ From Proposition \ref{Nil} we have that ${\rm NF}(x_{3}^{(i-j)a_{3}}-x_{2}^{(i-j)a_{32}}x_{4}^{(i-j)a_{34}}|T)$=0 and therefore ${\rm NF}({\rm spoly}(g_{i},h_{j})|T)=0$.
\end{enumerate}

We show that ${\rm NF}({\rm spoly}(g_{i},p)|T)=0$, for every $0 \leq i \leq d-2$. First assume that $a_{4}-da_{34}>0$. Then ${\rm spoly}(g_{i},p)=x_{1}^{a_{21}}x_{4}^{(i+1)a_{34}} \left[(x_{2}^{a_{32}}x_{4}^{a_{34}})^{d-i-1}-(x_{3}^{a_{3}})^{d-i-1} \right]$  From Proposition \ref{Nil} we have that ${\rm NF}(x_{2}^{(d-i-1)a_{32}}x_{4}^{(d-i-1)a_{34}}-x_{3}^{(d-i-1)a_{3}}|T)$=0 and therefore ${\rm NF}({\rm spoly}(g_{i},p)|T)=0$.

Suppose now that $a_{4}-da_{34} \leq 0$. Then ${\rm spoly}(g_{i},p)=x_{1}^{a_1}x_{2}^{a_{2}+(d-i-2)a_{32}}x_{4}^{da_{34}-a_{4}}-x_{1}^{a_{21}}x_{3}^{(d-i-1)a_{3}}x_{4}^{(i+1)a_{34}}$. Let $i=d-2$, then ${\rm spoly}(g_{d-2},p)=x_{1}^{a_{21}}x_{4}^{da_{34}-a_{4}}(x_{1}^{a_{41}}x_{2}^{a_2}-x_{3}^{a_3}x_{4}^{a_{14}})$. We already proved that ${\rm NF}(x_{1}^{a_{41}}x_{2}^{a_2}-x_{3}^{a_3}x_{4}^{a_{14}}|T)=0$, so ${\rm NF}({\rm spoly}(g_{d-2},p)|T)=0$. Suppose that $0 \leq i \leq d-3$. We distinguish the following cases: \begin{enumerate} \item ${\rm LM}({\rm spoly}(g_{i},p))=x_{1}^{a_1}x_{2}^{a_{2}+(d-i-2)a_{32}}x_{4}^{da_{34}-a_{4}}$. Notice that $da_{34}-a_{4}<a_{34}$ by definition of $d$. Only ${\rm LM}(f_2)$ divides ${\rm LM}({\rm spoly}(g_{i},p))$ and $${\rm spoly}(f_{2},{\rm spoly}(g_{i},p))=x_{1}^{a_{21}}x_{3}^{(d-i-1)a_{3}}x_{4}^{(i+1)a_{34}}-x_{1}^{a_{1}+a_{21}}x_{2}^{(d-i-2)a_{32}}x_{3}^{a_{23}}x_{4}^{da_{34}-a_{4}}.$$ Then ${\rm LM}({\rm spoly}(f_{2},{\rm spoly}(g_{i},p)))=x_{1}^{a_{21}}x_{3}^{(d-i-1)a_{3}}x_{4}^{(i+1)a_{34}},$ since $(d-i-1)(a_{3}-a_{32}-a_{34})+a_{4}+a_{32}-a_{1}-a_{23} \leq 0$. Notice that $a_{14}>(i+1)a_{34}$ for every $0 \leq i \leq d-3$, since $i+2<d$ and $a_{4}-(i+2)a_{34}>0$ by definition of $d$. Only ${\rm LM}(h_{d-i-2})$ divides ${\rm LM}({\rm spoly}(f_{2},{\rm spoly}(g_{i},p)))$ and also ${\rm ecart}({\rm spoly}(f_{2},{\rm spoly}(g_{i},p)))={\rm ecart}(h_{d-i-2})$. We have that $${\rm spoly}(h_{d-i-2},{\rm spoly}(f_{2},{\rm spoly}(g_{i},p)))=0$$ and ${\rm NF}({\rm spoly}(g_{i},p)|T)=0$.
\item ${\rm LM}({\rm spoly}(g_{i},p))=x_{1}^{a_{21}}x_{3}^{(d-i-1)a_{3}}x_{4}^{(i+1)a_{34}}$. Only the monomial ${\rm LM}(h_{d-i-2})$ divides ${\rm LM}({\rm spoly}(g_{i},p))$ and $${\rm spoly}(h_{d-i-2},{\rm spoly}(g_{i},p))=x_{1}^{a_1}x_{2}^{a_{2}+(d-i-2)a_{32}}x_{4}^{da_{34}-a_{4}}-x_{1}^{a_{1}+a_{21}}x_{3}^{(d-i-2)a_{32}}x_{3}^{a_{23}}x_{4}^{da_{34}-a_{4}}.$$ Then ${\rm LM}({\rm spoly}(h_{d-i-2},{\rm spoly}(g_{i},p)))=x_{1}^{a_1}x_{2}^{a_{2}+(d-i-2)a_{32}}x_{4}^{da_{34}-a_{4}}$, since $a_{2} \leq a_{21}+a_{23}$. Only ${\rm LM}(f_2)$ divides ${\rm LM}({\rm spoly}(h_{d-i-2},{\rm spoly}(g_{i},p)))$ and also ${\rm ecart}({\rm spoly}(h_{d-i-2},{\rm spoly}(g_{i},p)))={\rm ecart}(f_{2})$. We have that $${\rm spoly}(f_{2},{\rm spoly}(h_{d-i-2},{\rm spoly}(g_{i},p)))=0$$ and ${\rm NF}({\rm spoly}(g_{i},p)|T)=0$.
\end{enumerate}

We show that ${\rm NF}({\rm spoly}(h_{i},h_{j})|T)=0$, for every $i,j \in \{0,1,\ldots,d-2\}$ with $i<j$. We have that ${\rm spoly}(h_{i},h_{j})=x_{1}^{a_1}x_{2}^{ia_{32}} \left[(x_{2}^{a_{32}}x_{4}^{a_{34}})^{j-i}-(x_{3}^{a_3})^{j-i} \right]$. From Proposition \ref{Nil} we have that ${\rm NF}(x_{2}^{(j-i)a_{32}}x_{4}^{(j-i)a_{34}}-x_{3}^{(j-i)a_{3}}|T)$=0 and therefore ${\rm NF}({\rm spoly}(h_{i},h_{j})|T)=0$.

Finally we show that ${\rm NF}({\rm spoly}(h_{i},p)|T)=0$, for every $0 \leq i \leq d-2$. Suppose first that $a_{4}-da_{34} \leq 0$. Then ${\rm spoly}(h_{i},p)=x_{1}^{a_1}x_{2}^{ia_{32}} \left[(x_{2}^{a_{32}}x_{4}^{a_{34}})^{d-i-1}-(x_{3}^{a_3})^{d-i-1} \right]$. From Proposition \ref{Nil} we have that ${\rm NF}(x_{2}^{(d-i-1)a_{32}}x_{4}^{(d-i-1)a_{34}}-x_{3}^{(d-i-1)a_{3}}|T)$=0 and therefore ${\rm NF}({\rm spoly}(h_{i},p)|T)=0$. Suppose now that $a_{4}-da_{34}>0$, so $a_{2}-da_{32} \leq 0$ by definition of $d$. Then ${\rm spoly}(h_{i},p)=x_{1}^{a_{21}}x_{2}^{da_{32}-a_{2}}x_{4}^{a_{4}+(d-i-1)a_{34}}-x_{1}^{a_{1}}x_{2}^{ia_{32}}x_{3}^{(d-i-1)a_{3}}$. We have that ${\rm LM}({\rm spoly}(h_{i},p))=x_{1}^{a_{21}}x_{2}^{da_{32}-a_{2}}x_{4}^{a_{4}+(d-i-1)a_{34}}$, since $a_{4} \leq a_{41}+a_{42}$ and $a_{3}>a_{32}+a_{34}$. Notice that $da_{32}-a_{2}<a_{32}$ by definition of $d$. Only ${\rm LM}(f_4)$ divides ${\rm LM}({\rm spoly}(h_{i},p))$ and also ${\rm ecart}({\rm spoly}(h_{i},p))>{\rm ecart}(f_4)$. Then ${\rm spoly}(f_4,{\rm spoly}(h_i,p))=x_{1}^{a_1}x_{2}^{ia_{32}} \left[(x_{3}^{a_{3}})^{d-i-1}-(x_{2}^{a_{32}}x_{4}^{a_{34}})^{d-i-1} \right]$. From Proposition \ref{Nil} we have that ${\rm NF}(x_{3}^{(d-i-1)a_{3}}-x_{2}^{(d-i-1)a_{32}}x_{4}^{(d-i-1)a_{34}}|T)$=0 and therefore ${\rm NF}({\rm spoly}(h_{i},p)|T)=0$. \hfill $\square$\\

\begin{thm1} \label{CM2} Let $n_{1}={\rm min}\{n_{1},\ldots,n_{4}\}$ and $a_{3}>a_{32}+a_{34}$. Then ${\rm gr}_{\mathfrak{m}}(R)$ is Cohen-Macaulay if and only if the following conditions hold: \begin{enumerate}
\item $a_{2} \leq a_{21}+a_{23}$;
\item $a_{4} \leq a_{41}+a_{42}$;
\item $a_{42}+a_{13} \leq a_{21}+a_{34}$;
\item $(d-1)(a_{3}-a_{32}-a_{34})+a_{2}-a_{21}-a_{23} \leq 0$;
\item $(d-1)(a_{3}-a_{32}-a_{34})+a_{4}+a_{32}-a_{1}-a_{23} \leq 0$;
\item $d(a_{3}-a_{32}-a_{34})+a_{2}-a_{21}-a_{23} \leq 0$ when $a_{4}-da_{34}>0$;
\item $d(a_{3}-a_{32}-a_{34})+a_{4}+a_{32}-a_{1}-a_{23} \leq 0$ when $a_{4}-da_{34} \leq 0$.
\end{enumerate}
\end{thm1}

\noindent \textbf{Proof.} ($\Leftarrow$) Suppose that (1), (2), (3), (4), (5), (6) and (7) are true. Notice that $a_{1}>a_{13}+a_{14}$ since $n_{1}={\rm min}\{n_{1},\ldots,n_{4}\}$. By Proposition \ref{Basic2}, $T$ is a minimal standard basis for $I(C)$ with respect to the negative degree reverse lexicographic term ordering $<$ with $x_{4}>x_{3}>x_{2}>x_{1}$. Since $x_1$ does not divide the leading monomial of any polynomial in $T$, we have from Lemma \ref{AMS} that ${\rm gr}_{\mathfrak{m}}(R)$ is Cohen-Macaulay.\\
 ($\Rightarrow$) Suppose that ${\rm gr}_{\mathfrak{m}}(R)$ is Cohen-Macaulay. Let $S$ be a minimal standard basis of the ideal $I(C)$ with respect to a negative degree reverse lexicographical ordering $<$ that makes $x_1$ the lowest variable. The binomials $f_i$, $1 \leq i \leq 5$, are indispensable of $I(C)$, so they belong to $S$. By Lemma \ref{AMS}, $x_1$ does not divide ${\rm LM}(f_{i})$ for every $2 \leq i \leq 5$. Thus ${\rm LM}(f_{2})=x_{2}^{a_{2}}$, ${\rm LM}(f_{4})=x_{4}^{a_4}$ and ${\rm LM}(f_{5})=x_{2}^{a_{42}}x_3^{a_{13}}$. Therefore (1), (2) and (3) are true. Suppose that $(d-1)(a_{3}-a_{32}-a_{34})+a_{2}-a_{21}-a_{23}>0$. Then the binomial $g_{d-2}$ belongs to $I(C)$ and ${\rm LM}(g_{d-2})=x_{1}^{a_{21}}x_{4}^{(d-1)a_{34}}$. There is a binomial $B \in S$ such that ${\rm LM}(B)$ divides $x_{1}^{a_{21}}x_{4}^{(d-1)a_{34}}$. But ${\rm gr}_{\mathfrak{m}}(R)$ is Cohen-Macaulay, so $x_{1}$ does not divide ${\rm LM}(B)$. Thus ${\rm LM}(B)=x_{4}^{b_4}$ where $b_{4} \leq (d-1)a_{34}$. Notice that $(d-1)a_{34}<a_{4}$ by definition of $d$. So $b_{4}<a_{4}$ and therefore ${\rm LM}(B)$ strictly divides the monomial $x_{4}^{a_4}$, a contradiction to the fact that $x_{4}^{a_4}$ is an indispensable monomial. Suppose that $(d-1)(a_{3}-a_{32}-a_{34})+a_{4}+a_{32}-a_{1}-a_{23}>0$. Then the binomial $h_{d-2}$ belongs to $I(C)$ and ${\rm LM}(h_{d-2})=x_{1}^{a_{1}}x_{2}^{(d-2)a_{32}}$. There is a binomial $u \in S$ such that ${\rm LM}(u)$ divides $x_{1}^{a_{1}}x_{2}^{(d-2)a_{32}}$. But ${\rm gr}_{\mathfrak{m}}(R)$ is Cohen-Macaulay, so $x_{1}$ does not divide ${\rm LM}(u)$. Thus ${\rm LM}(u)=x_{2}^{c_2}$ where $c_{2} \leq (d-2)a_{32}$. Notice that $(d-2)a_{32}<a_{2}$ by definition of $d$. So $c_{2}<a_{2}$ and therefore ${\rm LM}(u)$ strictly divides the monomial $x_{2}^{a_2}$, a contradiction to the fact that $x_{2}^{a_2}$ is an indispensable monomial. Suppose that $a_{4}-da_{34}>0$ and let $d(a_{3}-a_{32}-a_{34})+a_{2}-a_{21}-a_{23}>0$. Then the binomial $p$ belongs to $I(C)$ and ${\rm LM}(p)=x_{1}^{a_{21}}x_{2}^{da_{32}-a_{2}}x_{4}^{da_{34}}$. Also there is a binomial $v \in S$ such that ${\rm LM}(v)$ divides $x_{1}^{a_{21}}x_{2}^{da_{32}-a_{2}}x_{4}^{da_{34}}$. Since ${\rm gr}_{\mathfrak{m}}(R)$ is Cohen-Macaulay, we have that ${\rm LM}(v)=x_{2}^{w_2}x_{4}^{w_4}$ where $w_2 \leq da_{32}-a_{2}$ and $w_4 \leq da_{34}$. Notice that $da_{32}-a_{2}<a_{32}$ by definition of $d$. So $w_{2}<a_{32}$ and $w_{4}<a_{4}$. Furthermore $w_{2}<a_{2}$ since $a_{32}<a_{2}$. But $v$ belongs to the ideal $I(C)$ and $\{f_{1},\ldots,f_{5}\}$ is a generating set for $I(C)$. Thus ${\rm LM}(v)$ is in the ideal generated by the monomials $x_3^{a_{13}} x_4^{a_{14}}, x_1^{a_1}, x_{2}^{a_2}, x_{1}^{a_{21}}x_{3}^{a_{23}}, x_{2}^{a_{32}}x_{4}^{a_{34}}, x_3^{a_{3}}, x_{4}^{a_4}, x_{1}^{a_{41}}x_{2}^{a_{42}}, x_{2}^{a_{42}}x_3^{a_{13}}$ and $x_{1}^{a_{21}}x_{4}^{a_{34}}$. So ${\rm LM}(v)$ is divided by at least one of the monomials $x_{2}^{a_2}$, $x_{2}^{a_{32}}x_{4}^{a_{34}}$ and $x_{4}^{a_4}$, a contradiction. Suppose that $a_{4}-da_{34} \leq 0$ and let $d(a_{3}-a_{32}-a_{34})+a_{4}+a_{32}-a_{1}-a_{23}>0$. Then the binomial $p$ belongs to $I(C)$ and ${\rm LM}(p)=x_{1}^{a_{1}}x_{2}^{(d-1)a_{32}}x_{4}^{da_{34}-a_{4}}$. There is a binomial $q \in S$ such that ${\rm LM}(q)$ divides $x_{1}^{a_{1}}x_{2}^{(d-1)a_{32}}x_{4}^{da_{34}-a_{4}}$. Since ${\rm gr}_{\mathfrak{m}}(R)$ is Cohen-Macaulay, we have that ${\rm LM}(q)=x_{2}^{e_2}x_{4}^{e_4}$ where $e_2 \leq (d-1)a_{32}$ and $e_4 \leq da_{34}-a_{4}$. Notice that $(d-1)a_{32}<a_{2}$ and $da_{34}-a_{4}<a_{34}$ by definition of $d$. So $e_{2}<a_{2}$ and $e_{4}<a_{34}$. Furthermore $e_{4}<a_{4}$ since $a_{34}<a_{4}$. But $q$ belongs to the ideal $I(C)$ and $\{f_{1},\ldots,f_{5}\}$ is a generating set for $I(C)$. Thus ${\rm LM}(q)$ is divided by at least one of the monomials $x_{2}^{a_2}$, $x_{2}^{a_{32}}x_{4}^{a_{34}}$ and $x_{4}^{a_4}$, a contradiction. \hfill $\square$\\

\begin{rem1} \label{basicrem2} {\rm Suppose that $I(C)$ is given as in case 2(b) of \cite{AKN} and also $a_{3}>a_{32}+a_{34}$. Note that $a_{4} \leq a_{41}+a_{43}$ since $n_{1}<n_{2}<n_{3}<n_{4}$. Then ${\rm gr}_{\mathfrak{m}}(R)$ is Cohen-Macaulay if and only if the following conditions hold: \begin{enumerate}
\item $a_{2} \leq a_{21}+a_{24}$;
\item $a_{24}+a_{13} \leq a_{41}+a_{32}$;
\item $(d-1)(a_{3}-a_{32}-a_{34})+a_{4}-a_{41}-a_{43} \leq 0$;
\item $(d-1)(a_{3}-a_{32}-a_{34})+a_{2}+a_{34}-a_{1}-a_{43} \leq 0$;
\item $d(a_{3}-a_{32}-a_{34})+a_{4}-a_{41}-a_{43} \leq 0$ when $a_{2}-da_{32}>0$;
\item $d(a_{3}-a_{32}-a_{34})+a_{2}+a_{34}-a_{1}-a_{43} \leq 0$ when $a_{2}-da_{32} \leq 0$.
\end{enumerate}}

\end{rem1}

\begin{rem1} \label{basicrem3} {\rm Suppose that $I(C)$ is given as in case 2(a) of \cite{AKN}. Note that $a_{2}>a_{23}+a_{24}$ and $a_{4} \leq a_{41}+a_{42}$ since $n_{1}<n_{2}<n_{3}<n_{4}$. Let $d$ be the minimum of $$\{i \in \mathbb{N}|a_{3}-ia_{23} \leq 0\} \cup \{i \in \mathbb{N}|a_{4}-ia_{24} \leq 0\}.$$ Then ${\rm gr}_{\mathfrak{m}}(R)$ is Cohen-Macaulay if and only if the following conditions hold: \begin{enumerate}
\item $a_{3} \leq a_{31}+a_{34}$;
\item $a_{34}+a_{12} \leq a_{41}+a_{23}$;
\item $(d-1)(a_{2}-a_{23}-a_{24})+a_{4}-a_{41}-a_{42} \leq 0$;
\item $(d-1)(a_{2}-a_{23}-a_{24})+a_{3}+a_{24}-a_{1}-a_{42} \leq 0$;
\item $d(a_{2}-a_{23}-a_{24})+a_{4}-a_{41}-a_{42} \leq 0$ when $a_{3}-da_{23}>0$;
\item $d(a_{2}-a_{23}-a_{24})+a_{3}+a_{24}-a_{1}-a_{42} \leq 0$ when $a_{3}-da_{23} \leq 0$.
\end{enumerate}}

\end{rem1}

\begin{rem1} \label{basicrem4} {\rm Suppose that $I(C)$ is given as in case 3(b) of \cite{AKN}. Note that $a_{2}>a_{23}+a_{24}$, $a_{3} \leq a_{31}+a_{32}$ and $a_{4} \leq a_{41}+a_{43}$ since $n_{1}<n_{2}<n_{3}<n_{4}$. Let $d$ be the minimum of $$\{i \in \mathbb{N}|a_{3}-ia_{23} \leq 0\} \cup \{i \in \mathbb{N}|a_{4}-ia_{24} \leq 0\}.$$ Then ${\rm gr}_{\mathfrak{m}}(R)$ is Cohen-Macaulay if and only if the following conditions hold: \begin{enumerate}
\item $a_{43}+a_{12} \leq a_{31}+a_{24}$;
\item $(d-1)(a_{2}-a_{23}-a_{24})+a_{3}-a_{31}-a_{32} \leq 0$;
\item $(d-1)(a_{2}-a_{23}-a_{24})+a_{4}+a_{23}-a_{1}-a_{32} \leq 0$;
\item $d(a_{2}-a_{23}-a_{24})+a_{3}-a_{31}-a_{32} \leq 0$ when $a_{4}-da_{24}>0$;
\item $d(a_{2}-a_{23}-a_{24})+a_{4}+a_{23}-a_{1}-a_{32} \leq 0$ when $a_{4}-da_{24} \leq 0$.
\end{enumerate}}

\end{rem1}

\begin{ex1} {\rm Consider $n_{1}=1673$, $n_{2}=2236$, $n_{3}=2248$, and $n_{4}=2828$. Then $C$ is a Gorenstein noncomplete intersection monomial curve and $I(C)$ is minimally generated by $x_{3}^{7}x_{4}^{11}-x_{1}^{28}, x_{2}^{13}-x_{1}^{12}x_{3}^{4}, x_{3}^{11}-x_{2}^{6}x_{4}^{4}, x_{4}^{15}-x_{1}^{16}x_{2}^{7}, x_{2}^{7}x_{3}^{7}-x_{1}^{12}x_{4}^4$. Here $d=3$ and $a_{4}-da_{34}=3>0$. We have that $d(a_{3}-a_{32}-a_{34})+a_{2}-a_{21}-a_{23}=0$ and $(d-1)(a_{3}-a_{32}-a_{34})+a_{4}+a_{32}-a_{1}-a_{23}=-9< 0$. By Theorem \ref{CM2}, the ring ${\rm gr}_{\mathfrak{m}}(R)$ is Cohen-Macaulay. Furthermore, from \cite[Lemma 5.5.11]{GP} we have that $T_{\ast}=\{x_{3}^{7}x_{4}^{11},x_{2}^{13},x_{2}^{6}x_{4}^{4},x_{4}^{15},x_{2}^{7}x_{3}^{7}, x_{2}x_{3}^{18},x_{3}^{18}x_{4}^{7}, x_{3}^{29}-x_{1}^{12}x_{2}^{5}x_{4}^{12}\}$ is a generating set for $I(C)_{\ast}$. }
\end{ex1}

\begin{ex1} {\rm Consider $n_{1}=890$, $n_{2}=1944$, $n_{3}=933$, and $n_{4}=1275$. Then $C$ is a Gorenstein noncomplete intersection monomial curve and $I(C)$ is minimally generated by $x_{3}^{5}x_{4}^{11}-x_{1}^{21}, x_{2}^{7}-x_{1}^{9}x_{3}^{6},x_{3}^{11}- x_{2}^{2}x_{4}^{5}, x_{4}^{16}-x_{1}^{12}x_{2}^{5}, x_{2}^{5}x_{3}^{5}-x_{1}^{9}x_{4}^5$. Here $d=4$ and $(d-1)(a_{3}-a_{32}-a_{34})+a_{2}-a_{21}-a_{23}=4>0$, so from Theorem \ref{CM2} the ring ${\rm gr}_{\mathfrak{m}}(R)$ is not Cohen-Macaulay.}
\end{ex1}

\begin{thm1} Let $n_{1}={\rm min}\{n_{1},\ldots,n_{4}\}$ and $C$ be a Gorenstein noncomplete intersection monomial curve such that ${\rm gr}_{\mathfrak{m}}(R)$ is Cohen-Macaulay. \begin{enumerate} \item If $a_{3} \leq a_{32}+a_{34}$, then $\mu(I(C)_{\ast})=5$. \item If $a_{3}>a_{32}+a_{34}$, then $\mu(I(C)_{\ast})=2d+2$.
\end{enumerate}
\end{thm1}
\noindent \textbf{Proof.}  Let $\pi: K[x_{1},x_{2},x_{3},x_{4}] \rightarrow K[y_1,x_{2},x_{3},x_{4}]$ be the $K$-algebra homomorphism defined by $\pi(x_{1})=y_{1}=0$ and $\pi(x_i)=x_{i}$, for $2 \leq i \leq 4$. From \cite[Lemma 1.2]{SS} we have that $\mu(I(C)_{\ast})$ equals the minimal number of generators of $\pi(I(C)_{\ast})$.\\
(1) By Proposition \ref{Basic1}, the set $G$ is a standard basis for $I(C)$, so, from \cite[Lemma 5.5.11]{GP}, $I(C)_{\ast}$ is generated by the least homogeneous summands of the elements in $G$. Thus $\pi(I(C)_{\ast})$ is minimally generated by either $\{x_3^{a_{13}} x_4^{a_{14}}, x_{2}^{a_2} ,x_{3}^{a_{3}}, x_{4}^{a_4}, x_{2}^{a_{42}}x_3^{a_{13}}\}$ or $\{x_3^{a_{13}} x_4^{a_{14}}, x_{2}^{a_2}, x_{3}^{a_{3}}-x_{2}^{a_{32}}x_{4}^{a_{34}}, x_{4}^{a_4}, x_{2}^{a_{42}}x_3^{a_{13}}\}$, and therefore $\mu(I(C)_{\ast})=5$.\\ (2) By Proposition \ref{Basic2}, the set $T$ is a standard basis for $I(C)$, so, from \cite[Lemma 5.5.11]{GP}, $I(C)_{\ast}$ is generated by the least homogeneous summands of the elements in $T$. Thus $\pi(I(C)_{\ast})$ is minimally generated by  $$\{x_{2}^{a_2}, x_{2}^{a_{32}}x_{4}^{a_{34}}, x_{4}^{a_4}\} \cup \{x_{2}^{a_{2}-(i+1)a_{32}}x_{3}^{ia_{3}+a_{13}}\}| 0 \leq i \leq d-2\} \cup$$ $$\{x_{3}^{ia_{3}+a_{13}}x_{4}^{a_{4}-(i+1)a_{34}}| 0 \leq i \leq d-2\} \cup \{x_{3}^{(d-1)a_{3}+a_{13}}\}$$ and therefore $\mu(I(C)_{\ast})=2d+2$. \hfill $\square$\\

The next proposition provides a family of Gorenstein noncomplete intersection monomial curves $C$ such that ${\rm gr}_{\mathfrak{m}}(R)$ is Cohen-Macaulay and $\mu(I(C)_{\ast})$ is large.

\begin{prop1} \label{infinite} Let $c \geq 2$ be an integer such that $227$ does not divide $c+80$, and let $n_{1}=3c^{2}-4c+2$, $n_{2}=6c^{2}+9c+2$, $n_{3}=4c^2+5c-3$, $n_{4}=6c^{2}+6c-11$. Then ${\rm gr}_{\mathfrak{m}}(R)$ is Cohen-Macaulay and $\mu(I(C)_{\ast})=2c+2$.
\end{prop1}

\noindent \textbf{Proof.} We have that $${\rm gcd}(n_{2},n_{4})={\rm gcd}(n_{2},n_{2}-n_{4})={\rm gcd}(n_{2}-2c(n_{2}-n_{4}),n_{2}-n_{4})=$$ $${\rm gcd}(n_{2}-2c(n_{2}-n_{4})+6(n_{2}-n_{4}),n_{2}-n_{4})={\rm gcd}(c+80,3c+13)={\rm gcd}(c+80,227)=1,$$ since $227$ is a prime integer which does not divide $c+80$. So ${\rm gcd}(n_{1},n_{2},n_{3},n_{4})=1$. By \cite[Theorem 4]{GS}, $C$ is a Gorenstein noncomplete intersection monomial curve and $I(C)$ is minimally generated by the set $$\{x_{3}x_{4}^{c-1}-x_{1}^{2c+4},x_{2}^{c}-x_{1}^{2c+3}x_{3}^{2},x_{3}^{3}-x_{2}x_{4},x_{4}^{c}-x_{1}x_{2}^{c-1},x_{2}^{c-1}x_{3}-x_{1}^{2c+3}x_{4}\}.$$ Here $d=c$ and $a_{4}-da_{34}=0$. We have that $(d-1)(a_{3}-a_{32}-a_{34})+a_{2}-a_{21}-a_{23}=-6<0$ and also $d(a_{3}-a_{32}-a_{34})+a_{4}+a_{32}-a_{1}-a_{23}=-5<0$. By Theorem \ref{CM2}, ${\rm gr}_{\mathfrak{m}}(R)$ is Cohen-Macaulay. From Proposition \ref{Basic2} we have that the set $$T=\{x_{2}^{c}-x_{1}^{2c+3}x_{3}^{2}, x_{3}^{3}-x_{2}x_{4}, x_{4}^{c}-x_{1}x_{2}^{c-1}\} \cup \{x_{2}^{c-(i+1)}x_{3}^{3i+1}-x_{1}^{2c+3}x_{4}^{i+1}| 0 \leq i \leq c-2\} \cup$$ $$\{x_{3}^{3i+1}x_{4}^{c-(i+1)}-x_{1}^{2c+4}x_{2}^{i}| 0 \leq i \leq c-2\} \cup \{x_{3}^{3c-2}-x_{1}^{2c+4}x_{2}^{c-1}\}$$ is a standard basis for $I(C)$ with respect to the negative degree reverse lexicographic term ordering $<$ with $x_{4}>x_{3}>x_{2}>x_{1}$. From \cite[Lemma 5.5.11]{GP} the ideal $I(C)_{\ast}$ is generated by the least homogeneous summands of the elements in $T$. Then $I(C)_{\ast}$ is minimally generated by the set $$\{x_{2}^{c}, x_{2}x_{4}, x_{4}^{c}-x_{1}x_{2}^{c-1}\} \cup \{x_{2}^{c-(i+1)}x_{3}^{3i+1}| 0 \leq i \leq c-2\} \cup$$ $$ \{x_{3}^{3i+1}x_{4}^{c-(i+1)}| 0 \leq i \leq c-2\} \cup \{x_{3}^{3c-2}\}.$$ Thus $\mu(I(C)_{\ast})=2c+2$. \hfill $\square$\\

\textbf{Acknowledgments.} The author would like to thank the referee for careful reading and helpful comments.


\begin{thebibliography}{50}

\bibitem{AKN} F. Arslan, A. Katsabekis, M. Nalbandiyan. {\em On the Cohen-Macaulayness of tangent cones of monomial curves in $A^{4}(K)$}, Turk. J. Math. {\bf 43} (2014), 1425-1446.
\bibitem{AM} F. Arslan, P. Mete, {\em Hilbert functions of Gorenstein monomial curves}, Proc. Amer. Math. Soc. {\bf 135} (2007), 1993--2002.

\bibitem{AMS} F. Arslan, P. Mete, M. Sahin. {\em Gluing and Hilbert functions of monomial curves}, Proc. Amer. Math. Soc. {\bf 137} (2009), 2225–2232. 
	
\bibitem{Bresinsky75} H. Bresinsky. {\em Symmetric semigroups of integers generated by 4 elements}, Manuscripta Math. {\bf 17} (3) (1975), 205–219. 
\bibitem{GS} P. Gimenez, H. Srinivasan. {\em A note on Gorenstein monomial curves},  Bull. Braz. Math. Soc. (N.S.) {\bf 45} (2014), no. 4, 671–678.
\bibitem{GP} G. M. Greuel, G. Pfister, {\em A Singular Introduction to Commutative Algebra}, Springer-Verlag, 2002.
\bibitem{HerSta} J. Herzog, D. Stamate. {\em On the defining equations of the tangent cone of a numerical semigroup ring}, J. Algebra {\bf 418} (2014), 8–28.
\bibitem{katsabekis} A. Katsabekis, {\em Hilbert series of tangent cones for Gorenstein monomial curves in $A^{4}(K)$}, Turk. J. Math. {\bf 45} (2021), 597--616. 
\bibitem{KO} A. Katsabekis, I. Ojeda. {\em An indispensable classification of monomial curves in $\mathbb{A}^{4}(\mathbb{K})$}, Pacific J. Math. {\bf 268} (2014), 95–116. 
\bibitem{LTV} N. P. H. Lan, N. C. Tu, T. Vu. {\em Betti numbers of the tangent cones of monomial space curves}, arXiv: 2307.05589.
\bibitem{SS} M. Sahin, N. Sahin, {\em Betti numbers for certain Cohen-Macaulay tangent cones}, Bulletin of the Australian Mathematical Soc. {\bf 99} (2019), 68-77.
\bibitem{Sahin} N. Sahin, {\em On the monotonicity of the Hilbert functions for 4-generated pseudo-symmetric monomial curves }, arXiv:2302.09356.
\bibitem{Shibuta} T. Shibuta, {\em Cohen-Macaulayness of almost complete intersection tangent cones}, J. Algebra {\bf 319} (8) (2008), 3222-3243.
\bibitem{Sta} R. P. Stanley, {\em Hilbert Functions of Graded Algebras}, Advances in Math. {\bf 28} (1978), 57-83.
\bibitem{Sturmfels95} B. Sturmfels. {\em Gr\"obner Bases and Convex Polytopes}, University Lecture Series Vol. 8. Providence, RI, USA: American Mathematical Society, 1996.
	
\end{thebibliography}
\end{document}